\documentclass[12pt]{article}

\begin{document}

\begin{titlepage}
\title{\bf Lifts of Lorentzian r- Paracontact Structure: \\A Geometrical Dynamics Meaning}
\author{ Mehmet Tekkoyun \footnote{tekkoyun@pamukkale.edu.tr} \\
 {\small Department of Mathematics, Pamukkale University,}\\
{\small 20070 Denizli, Turkey}}
\date{\today}
\maketitle

\begin{abstract}

The goal of this paper, using lifting theory it is to produce
almost paracomplex structures on the tangent bundle  of almost
Lorentzian r-paracontact manifold endowed with almost Lorentzian
r-paracontact structure. Finally, we discuss the effect  over
dynamics systems of the produced geometrical structures.



\end{abstract}
\end{titlepage}

\section{Introduction}

It is well known that the method of lift has an important role in modern
differentiable geometry. Because, differentiable structures and mechanical
systems defined on a manifold $M$ or a physical space can be lifted to the
same type of structures on its tangent and cotangent bundles which are
phase-spaces of velocities and momentums of a given configuration space. It
was studied lifting theory of complex and paracomplex structures in [1-4]
and the differential geometry of tangent and cotangent bundles in [5-7].

The paper is structured as follows. In section 2, we recall the fundamental
structures about paracomplex, paracontact, r-paracontact, Lorentzian
paracontact and Lorentzian r-paracontact manifolds given in [8-12]. In
section 3, we give vertical, complete and horizontal lifts of paracomplex
structures[4]. In section 4, we produce almost paracomplex structures by
means of complete and horizontal lifts of almost Lorentzian r-paracontact
structure on almost Lorentzian r-paracontact manifold. Also we shall
determine the influence on field of dynamics systems of differential
geometric elements.

Along this paper, all mappings and manifolds will be understood to be of
class differentiable and the sum is taken over repeated indices. In this
study we denote by $\Im _{s}^{r}(M)$ the set of all tensor fields of class $%
C^{\infty }$ and of type $(r,s)$ in $M$. We now put $\Im
(M)=\sum_{r,s=1}^{\infty }\Im _{s}^{r}(M)$; which is the set of all tensor
fields in $M$. For example $\Im _{0}^{0}(M),\Im _{0}^{1}(M),\Im _{1}^{0}(M)$
and $\Im _{1}^{1}(M)$ are the set of functions, vector fields,1-forms and
tensor fields of type (1,1) on $M$, respectively. We give by $\mathcal{F}(M)$
to $\Im _{0}^{0}(M),$ by $\chi (M)$ to $\Im _{0}^{1}(M)$ and by $\chi
^{*}(M) $ to $\Im _{1}^{0}(M).$ Similarly, we respectively denote by $\Im
_{s}^{r}(T(M))$ and $\Im (T(M))$ the corresponding sets of tensor fields in
the tangent bundle $T(M)$. Also $v$, $c$ and $h$ will denote the vertical,
complete and horizontal lifts to $TM$ of geometric structures on $M$,
respectively.

\section{\textbf{Basic Structures}}

\subsection{Paracomplex Geometry}

A tensor field of type (1,1) $J$ on $M$ such that $J^{2}=I$ is an \textit{%
almost product structure} $J$ on manifold $M$ of dimension $2m$. Hence, we
said to be an \textit{almost product manifold} the pair $(M,J)$. If the two
eigenbundles $T^{+}M$ and $T^{-}M$ associated to the eigenvalues +1 and -1
of $J$, respectively, have the same rank, then an almost product manifold $%
(M,J)$ is said an \textit{almost paracomplex manifold}. The dimension of an
almost paracomplex manifold is necessarily even. Equivalently, a splitting
of the tangent bundle $TM$ of manifold $M$, into the Whitney sum of two
subbundles on $T^{\pm }M$ of the same fiber dimension is called an \textit{%
almost paracomplex structure} on $M.$ An almost paracomplex structure on
manifold $M $ may alternatively be defined as a $G$- structure on $M$ with
structural group $GL(n,\mathbf{R})\times GL(n,\mathbf{R})$.

A \textit{paracomplex manifold is }an almost paracomplex manifold $(M,J)$
such that the $G$- structure defined by the tensor field $J$ is integrable. $%
(\ x^{\alpha },\,\ y^{\alpha }),1\leq \alpha \leq m$ be a real coordinate
system on a neighborhood $U$ of any point $p$ of $M,$ and $\{(\frac{\partial
}{\partial x^{\alpha }})_{p},(\frac{\partial }{\partial y^{\alpha }})_{p}\}$
and $\{(dx^{\alpha })_{p},(dy^{\alpha })_{p}\}$ natural bases over $\mathbf{R%
}$ of the tangent space $T_{p}M$ and the cotangent space $T_{p}^{*}M$ of $M,$
respectively$.$ Hence we shall define as:
\begin{equation}
J(\frac{\partial }{\partial x^{\alpha }})=\frac{\partial }{\partial
y^{\alpha }},\,J(\frac{\partial }{\partial y^{\alpha }})=\frac{\partial }{%
\partial x^{\alpha }} ,  \label{1.1}
\end{equation}
\begin{equation}
J^{*}(dx^{\alpha })=-dy^{\alpha },\,J^{*}(dy^{\alpha })=-dx^{\alpha }.
\label{1.2}
\end{equation}
Let $\ z^{\alpha }=x^{\alpha }+$\textbf{j}$y^{\alpha },\,\overline{z}%
^{\alpha }=x^{\alpha }-$\textbf{j}$y^{\alpha },\,1\leq \alpha \leq m\,,\,$%
\textbf{j}$^{2}=1,$ be a paracomplex local coordinate system on a
neighborhood $U$ of any point $p$ of $M.$ We define the vector fields and
the dual covector fields as:
\begin{equation}
(\frac{\partial }{\partial z^{\alpha }})_{p}=\frac{1}{2}\{(\frac{\partial }{%
\partial x^{\alpha }})_{p}-\mathbf{j}(\frac{\partial }{\partial y^{\alpha }}%
)_{p}\},\,(\frac{\partial }{\partial \overline{z}^{\alpha }})_{p}=\frac{1}{2}%
\{(\frac{\partial }{\partial x^{\alpha }})_{p}+\mathbf{j}(\frac{\partial }{%
\partial y^{\alpha }})_{p}\} ,  \label{1.3}
\end{equation}
\begin{equation}
\left( dz^{\alpha }\right) _{p}=\left( dx^{\alpha }\right) _{p}+\mathbf{j}%
(dy^{\alpha })_{p},\,\left( d\overline{z}^{\alpha }\right) _{p}=\left(
dx^{\alpha }\right) _{p}-\mathbf{j}(dy^{\alpha })_{p} ,  \label{1.4}
\end{equation}
which represent the bases of the tangent space $T_{p}M$ and cotangent space $%
T_{p}^{*}M$ of $M$, respectively. Then, using \textbf{j}$^{2}=1$ it is
obtained
\begin{equation}
J(\frac{\partial }{\partial z^{\alpha }})=-\mathbf{j}\frac{\partial }{%
\partial z^{\alpha }},\,J(\frac{\partial }{\partial \overline{z}^{\alpha }})=%
\mathbf{j}\frac{\partial }{\partial \overline{z}^{\alpha }},  \label{1.5}
\end{equation}
\begin{equation}
J^{*}(dz^{\alpha })=-\mathbf{j}dz^{\alpha },\,J^{*}(d\overline{z}^{\alpha })=%
\mathbf{j}d\overline{z}^{\alpha }.  \label{1.6}
\end{equation}
For each $\,p\in M,$ let $T_{p}M$ \thinspace be set of tangent vectors $%
Z_{p}=Z^{\alpha }(\frac{\partial }{\partial z^{\alpha }})_{p}+\overline{Z}%
^{\alpha }(\frac{\partial }{\partial \overline{z}^{\alpha }})_{p},$ then $%
TM\,$ is the union of these vector spaces. Thus, tangent bundle of a
paracomplex manifold $M$ is $\left( TM,\tau _{M},M\right) ,$ where canonical
projection $\tau _{M}$ is $\tau _{M}:TM\rightarrow M$ $(\,\tau
_{M}(\,Z_{p})=p)$ and, in addition this map are surjective submersion. Then
now, coordinates $\left\{ z_{{}}^{\alpha },\overline{z}_{{}}^{\alpha
},z_{{}}^{\acute{}\alpha },\overline{z}_{{}}^{\acute{}\alpha }\right\} ,$ $%
1\leq \alpha \leq m,$ are taken into considering as local coordinates for $%
TM $.

\subsection{Almost Paracontact (r-Paracontact) Geometry}

Let $\overline{M}$ be an $n$-dimensional differentiable manifold. If there
exist on $\overline{M}$ a (1,1)-tensor field $\varphi $, a vector field $\xi
$ and a 1-form $\eta $ satisfying

\begin{equation}
\eta (\xi )=1,\,\,\,\varphi ^{2}=I-\eta \otimes \xi ,  \label{1.7}
\end{equation}

where $I$ is the identity, then $\overline{M}$ is said to be an almost
paracontact manifold. In the almost paracontact manifold, the following
relations hold good:

\begin{equation}
\,\varphi (\xi )=0,\,\eta \circ \varphi =0,\,\,rank(\varphi )=n-1.
\label{1.8}
\end{equation}

Every almost paracontact manifold has a positive definite Riemannian metric $%
G$ such that

\begin{equation}
\eta (\overline{X})=G(\xi ,\overline{X}),\,\,G(\varphi \overline{X},\varphi
\overline{Y})=G(\overline{X},\overline{Y})-\eta (\overline{X})\eta (%
\overline{Y}),\,\,\overline{X},\overline{Y}\in \chi (\overline{M})
\label{1.9}
\end{equation}

where $\chi (\overline{M})$ denotes the set of differentiable vector fields
on $\overline{M}$. In this case, we say that $\overline{M}$ has an almost
paracontact Riemannian structure $(\varphi ,\xi ,\eta ,G)$ and $\overline{M}$
is said to be an almost paracontact Riemannian manifold.

This structure generalizes as follows. One takes care of a tensor field $%
\mathbf{F}$ of type (1,1) on a manifold $\overline{M}$ of dimension $(2n+r)$%
. If there exists on $\overline{M}$ the vector fields $(\xi _{\alpha })$ and
the 1- forms $(\eta ^{\alpha })$ such that

\begin{equation}
\eta ^{\alpha }(\xi _{\beta })=\delta _{\beta }^{\alpha },\,\,\mathbf{F}(\xi
_{\alpha })=0,\,\eta ^{\alpha }\circ {\mathbf{F}}=0,\,\,\,{\mathbf{F}}%
^{2}=I+\sum_{\alpha =1}^{r}\xi _{\alpha }\otimes \eta ^{\alpha },
\label{1.10}
\end{equation}

then the structure $(\mathbf{F},\xi _{\alpha },\eta ^{\alpha })$ is also
called a framed ${\mathbf{F}}$ $(3,-1)$-structure, where $(\alpha ,\beta
=1,2,...,r)$ and $\delta _{\beta }^{\alpha }$ denotes Kronecker delta. A
manifold $\overline{M}$ endowed with a framed ${\mathbf{F}}(3,-1)$-structure
is called an almost $r$-paracontact manifold. In addition, we say that the
slit tangent bundle of a Finsler space carries a natural framed ${\mathbf{F}}%
(3;-1)$- structure and then the framed $\mathbf{F}(3;-1)$- structure induces
an almost paracontact structure.

\subsection{Lorentzian Almost Paracontact (r-Paracontact) Geometry}

Let $\overline{M}$ be an $n$-dimensional differentiable manifold equipped
with a triple $(\varphi ,\xi ,\eta ),$ where $\varphi $ is a (1,1)-tensor
field$,\xi $ a vector field and $\eta $ is a 1-form on $\overline{M}$
satisfying

\begin{equation}
\eta (\xi )=-1,\,\,\,\varphi ^{2}=I+\eta \otimes \xi ,  \label{1.11}
\end{equation}

where $I$ is the identity. In manifold $\overline{M},$ the following
relations hold good:

\begin{equation}
\,\varphi (\xi )=0,\,\eta \circ \varphi =0,\,\,rank(\varphi )=n-1.
\label{1.12}
\end{equation}

If $\overline{M}$ has a Lorentzian metric $G$ such that

\begin{equation}
G(\varphi \overline{X},\varphi \overline{Y})=G(\overline{X},\overline{Y}%
)+\eta (\overline{X})\eta (\overline{Y}),\,\,\overline{X},\overline{Y}\in
\chi (\overline{M})  \label{1.13}
\end{equation}

we say that $\overline{M}$ has a Lorentzian almost paracontact structure $%
(\varphi ,\xi ,\eta ,G)$ and $\overline{M}$ is said to be a Lorentzian
almost paracontact manifold.

This structure generalizes as follows. One takes care of a tensor field $%
\mathbf{F}$ of type (1,1) on a manifold $\overline{M}$ of dimension $(2n+r)$%
. If there exists on $\overline{M}$ the vector fields $(\xi _{\alpha })$ and
the 1- forms $(\eta ^{\alpha })$ such that

\begin{equation}
\eta ^{\alpha }(\xi _{\beta })=-\delta _{\beta }^{\alpha },\,\,\mathbf{F}%
(\xi _{\alpha })=0,\,\eta ^{\alpha }\circ {\mathbf{F}}=0,\,\,\,{\mathbf{F}}%
^{2}=I-\sum_{i=1}^{r}\xi _{\alpha }\otimes \eta ^{\alpha },  \label{1.14}
\end{equation}

then the structure $(\mathbf{F},\xi _{\alpha },\eta ^{\alpha })$ is also
called a Lorentzian almost $r$-paracontact structure, where $(\alpha ,\beta
=1,2,...,r)$ and $\delta _{\beta }^{\alpha }$ denotes Kronecker delta. A
manifold $\overline{M}$ endowed with $(\mathbf{F},\xi _{\alpha },\eta
^{\alpha })$-structure is called a Lorentzian almost $r$-paracontact
manifold.

\section{\textbf{Lifting theory of Paracomplex Structures}}

\subsection{Vertical Lifts}

The \textit{vertical\thinspace \thinspace lift} of\thinspace \thinspace
paracomplex function $f$ $\in \mathcal{F}(M)$\thinspace to $TM$ is the
function $f^{v}$ in $\mathcal{F}(TM)$ given by
\begin{equation}
f^{v}=f\circ \tau _{M},  \label{2.1}
\end{equation}

where $\tau _{M}:TM\rightarrow M$ canonical projection.$\,\,$\thinspace We
have $rang(f^{v})=rang(f),$ since$\,\,$
\begin{equation}
f^{v}(Z_{p})=f(\tau _{M}(Z_{p}))=f(p),\,\,\,\,\,\,\,\forall \,Z_{p}\in TM.
\label{2.2}
\end{equation}

The \textit{vertical lift} of a vector field\textbf{\ }$Z\in \chi (M)$ to $%
TM $ is the vector field $Z^{v}\in \chi (TM)$ given by
\begin{equation}
Z^{v}(f^{c})=(Zf)^{v},\,\,\,\,\,\,\forall f\in \mathcal{F}(M).  \label{2.3}
\end{equation}

If $Z=Z^{\alpha }\frac{\partial }{\partial z^{\alpha }}+\overline{Z}^{\alpha
}\frac{\partial }{\partial \overline{z}^{\alpha }}$ we have
\begin{equation}
Z^{v}=(Z^{\alpha })^{v}\frac{\partial }{\partial z^{\acute{}\alpha }}+(%
\overline{Z}^{\alpha })^{v}\frac{\partial }{\partial \overline{z}^{\acute{}%
\alpha }},1\leq \alpha \leq m.  \label{2.4}
\end{equation}

The \textit{vertical lift} of a 1-form\textbf{\ }$\omega \in \chi ^{*}(M)$
to $TM$ is the 1-form $\omega ^{v}\in \chi ^{*}(TM)$ given by
\begin{equation}
\omega ^{v}(Z^{c})=(\omega Z)^{v},{\,\,\,\,\,\,}\forall Z\in \chi (M).
\label{2.5}
\end{equation}

The vertical lift of the paracomplex 1-form $\omega $ given by
\begin{equation}
\omega =\omega _{\alpha }dz^{\alpha }+\overline{\omega }_{\alpha }d\overline{%
z}^{\alpha }  \label{2.6}
\end{equation}
is
\begin{equation}
\omega ^{v}=(\omega _{\alpha })^{v}dz^{\alpha }+(\overline{\omega }_{\alpha
})^{v}d\overline{z}^{\alpha },1\leq \alpha \leq m.  \label{2.7}
\end{equation}

The \textit{vertical lift} of a paracomplex tensor field of type (1,1)%
\textbf{\ }$F\in \Im _{1}^{1}(M)$ to tangent bundle $TM$ of a paracomplex
manifold $M$ is the tensor field $F^{v}\in \Im _{1}^{1}(TM)$ given by
\begin{equation}
F^{v}(Z^{c})=(FZ)^{v},{\,\,\,\,\,\,}\forall Z\in \chi (M).  \label{2.8}
\end{equation}

Clearly, we have
\begin{equation}
F^{v}=(F_{\alpha }^{\beta })^{v}\frac{\partial }{\partial z^{\acute{}\beta }}%
\otimes dz^{\alpha }+(\overline{F}_{\alpha }^{\beta })^{v}\frac{\partial }{%
\partial \overline{z}^{\acute{}\beta }}\otimes d\overline{z}^{\alpha },1\leq
\alpha ,\beta \leq m,  \label{2.9}
\end{equation}

where $F=F_{\alpha }^{\beta }\frac{\partial }{\partial z^{\beta }}\otimes
dz^{\alpha }+\overline{F}_{\alpha }^{\beta }\frac{\partial }{\partial
\overline{z}^{\beta }}\otimes d\overline{z}^{\alpha }.$

The vertical lifts of paracomplex tensor fields have the general properties

\[
\begin{array}{ll}
i) & \,(f.g)^{v}=f^{v}.g^{v},(f+g)^{v}=f^{v}+g^{v},\, \\
ii) & (X+Y)^{v}=X^{v}+Y^{v},\,(fX)^{v}=f^{v}X^{v},\,X^{v}(f^{v})=0,\,\left[
X^{v},Y^{v}\right] =0, \\
iii) &
\begin{array}{l}
\chi (U)=Sp\left\{ \frac{\partial }{\partial z^{\alpha }},\frac{\partial }{%
\partial \overline{z}^{\alpha }}\right\} ,\,\chi (TU)=Sp\left\{ \frac{%
\partial }{\partial z^{\alpha }},\frac{\partial }{\partial \overline{z}%
^{\alpha }},\frac{\partial }{\partial z^{\acute{}\alpha }},\frac{\partial }{%
\partial \overline{z}^{\acute{}\alpha }}\right\} , \\
\;(\frac{\partial }{\partial z^{\alpha }})^{v}=\frac{\partial }{\partial z^{%
\acute{}\alpha }},(\frac{\partial }{\partial \overline{z}^{\alpha }})^{v}=%
\frac{\partial }{\partial \overline{z}^{\acute{}\alpha }},%
\end{array}
\\
iv) & \;(\omega +\theta )^{v}=\omega ^{v}+\theta ^{v},(f\omega
)^{v}=f^{v}\omega ^{v},(f\omega )^{v}=f^{v}\omega ^{v},\omega ^{v}(Z^{v})=0,
\\
v) &
\begin{array}{l}
\chi ^{*}(U)=Sp\left\{ dz^{\alpha },d\overline{z}^{\alpha }\right\} ,\chi
^{*}(TU)=Sp\left\{ dz^{\alpha },d\overline{z}^{\alpha },dz^{\acute{}\alpha
},d\overline{z}^{\acute{}\alpha }\right\} , \\
(dz^{\alpha })^{v}=\overline{d}z^{\alpha },(d\overline{z}^{\alpha })^{v}=%
\overline{d}\overline{z}^{\alpha },%
\end{array}
\\
vi) & \,\,F^{v}(Z^{v})=0,%
\end{array}
\]

for all\thinspace \thinspace \thinspace $f\in \mathcal{F}(M),$ $X,Y\in \chi
(M),\,\omega ,\theta \in \chi ^{*}(M),\,F\in \Im _{1}^{1}(M).$ Where $1\leq
\alpha \leq m$,$\,\,\left[ ,\right] $ is Lie bracket, $\overline{d}\,\,$%
denotes the differential operator on$\,\,TM$.

\subsection{Complete Lifts}

The \textit{complete \thinspace lift} of\thinspace \thinspace \thinspace
paracomplex function $f$ \thinspace $\in \mathcal{F}(M)$\ to $TM$ is the
function $f^{c}\in \mathcal{F}(TM)$ given by
\begin{equation}
f^{c}=z^{\acute{}\alpha }(\frac{\partial f}{\partial z^{\alpha }})^{v}+%
\overline{z}^{\acute{}\alpha }(\frac{\partial f}{\partial \overline{z}%
^{\alpha }})^{v},  \label{2.10}
\end{equation}
where we denote by $(z_{{}}^{\alpha },\overline{z}_{{}}^{\alpha },z_{{}}^{%
\acute{}\alpha },\overline{z}_{{}}^{\acute{}\alpha })$ local coordinates of
a chart-domain $TU\subset TM.$ Furthermore, for $\,Z_{p}\in TM$ we have
\begin{equation}
f^{c}(Z_{p})=z^{\acute{}\alpha }(Z_{p})(\frac{\partial f}{\partial z^{\alpha
}})^{v}(p)+\overline{z}^{\acute{}\alpha }(Z_{p})(\frac{\partial f}{\partial
\overline{z}^{\alpha }})^{v}(p).  \label{2.11}
\end{equation}

The \textit{complete lift} of a vector field\textbf{\ }$Z\in \chi (M)$ to $%
TM $ is the vector field $Z^{c}\in \chi (TM)$ given by
\begin{equation}
Z^{c}(f^{c})=(Zf)^{c},{\,\,\,\,\,\,}\forall f\in \mathcal{F}(M).
\label{2.12}
\end{equation}
Obviously, we obtain
\begin{equation}
Z^{c}=(Z^{\alpha })^{v}\frac{\partial }{\partial z^{\alpha }}+(\overline{Z}%
^{\alpha })^{v}\frac{\partial }{\partial \overline{z}^{\alpha }}+(Z^{\alpha
})^{c}\frac{\partial }{\partial z^{\acute{}\alpha }}+(\overline{Z}^{\alpha
})^{c}\frac{\partial }{\partial \overline{z}^{\acute{}\alpha }},1\leq \alpha
\leq m,  \label{2.13}
\end{equation}

where $Z=Z^{\alpha }\frac{\partial }{\partial z^{\alpha }}+\overline{Z}%
^{\alpha }\frac{\partial }{\partial \overline{z}^{\alpha }}.$

The \textit{complete lift} of a 1-form\textbf{\ }$\omega \in \chi ^{*}(M)$
to $TM$ is the 1-form $\omega ^{c}\in \chi ^{*}(TM)$ given by
\begin{equation}
\omega ^{c}(Z^{c})=(\omega Z)^{c},{\,\,\,\,\,\,}\forall Z\in \chi (M).
\label{2.14}
\end{equation}

If $\omega =\omega _{\alpha }dz^{\alpha }+\overline{\omega }_{\alpha }d%
\overline{z}^{\alpha }$ we calculate
\begin{equation}
\omega ^{c}=(\omega _{\alpha })^{c}dz^{\alpha }+(\overline{\omega }_{\alpha
})^{c}d\overline{z}^{\alpha }+(\omega _{\alpha })^{v}dz^{\acute{}\alpha }+(%
\overline{\omega }_{\alpha })^{v}d\overline{z}^{\acute{}\alpha },1\leq
\alpha \leq m.  \label{2.15}
\end{equation}

The \textit{complete lift} of a paracomplex tensor field of type (1,1)%
\textbf{\ }$F\in \Im _{1}^{1}(M)$ to $TM$ is the tensor field $F^{c}\in \Im
_{1}^{1}(TM)$ given by
\begin{equation}
F^{c}(Z^{c})=(FZ)^{c},{\,\,\,\,\,\,}\forall Z\in \chi (M).  \label{2.16}
\end{equation}

The complete lift of the paracomplex tensor field of type (1,1)\textbf{\ }$F$
is
\begin{equation}
\begin{array}{l}
F^{c}=(F_{\alpha }^{\beta })^{v}\frac{\partial }{\partial z^{\beta }}\otimes
dz^{\alpha }+(F_{\alpha }^{\beta })^{c}\frac{\partial }{\partial z^{\acute{}%
\beta }}\otimes dz^{\alpha }+(F_{\alpha }^{\beta })^{v}\frac{\partial }{%
\partial z^{\acute{}\beta }}\otimes dz^{\acute{\alpha}} \\
\,\,\,\,\,\,\,\,\,\,\,+(\overline{F}_{\alpha }^{\beta })^{v}\frac{\partial }{%
\partial \overline{z}^{\beta }}\otimes d\overline{z}^{\alpha }+(\overline{F}%
_{\alpha }^{\beta })^{c}\frac{\partial }{\partial \overline{z}^{\acute{}%
\beta }}\otimes d\overline{z}^{\alpha }+(\overline{F}_{\alpha }^{\beta })^{v}%
\frac{\partial }{\partial \overline{z}^{\acute{}\beta }}\otimes d\overline{z}%
^{\acute{}\alpha },%
\end{array}
\label{2.17}
\end{equation}

where $F=F_{\alpha }^{\beta }\frac{\partial }{\partial z^{\beta }}\otimes
dz^{\alpha }+\overline{F}_{\alpha }^{\beta }\frac{\partial }{\partial
\overline{z}^{\beta }}\otimes d\overline{z}^{\alpha }.$

The \textit{complete lift} of $J$ given by $J=\mathbf{j}\frac{\partial }{%
\partial z^{\alpha }}\otimes dz^{\alpha }-\mathbf{j}\frac{\partial }{%
\partial \overline{z}^{\alpha }}\otimes d\overline{z}^{\alpha }$ and being a
paracomplex tensor field of type (1,1) is

\begin{equation}
J^{c}=\mathbf{j}\frac{\partial }{\partial z^{\alpha }}\otimes dz^{\alpha }+%
\mathbf{j}\frac{\partial }{\partial z^{\acute{}\alpha }}\otimes dz^{\acute{}%
\alpha }-\mathbf{j}\frac{\partial }{\partial \overline{z}^{\alpha }}\otimes d%
\overline{z}^{\alpha }-\mathbf{j}\frac{\partial }{\partial \overline{z}^{%
\acute{}\alpha }}\otimes d\overline{z}^{\acute{}\alpha }.  \label{2.18}
\end{equation}

Because of $(J^{c})^{2}=I,$ $J^{c}$ is an almost paracomplex structure for
tangent bundle $TM$.

The complete lifts of paracomplex tensor fields obey the general properties

\[
\begin{array}{ll}
i) & (f+g)^{c}=f^{c}+g^{c},(f.g)^{c}=f^{c}.g^{v}+f^{v}.g^{c},\, \\
ii) &
\begin{array}{l}
\,(X+Y)^{c}=X^{c}+Y^{c},(fX)^{c}=f^{c}X^{v}+f^{v}X^{c}, \\
X^{c}(f^{v})=X^{v}(f^{c})=(Xf)^{v},X^{c}(f^{c})=(Xf)^{c}, \\
\left[ X^{v},Y^{c}\right] =\left[ X^{c},Y^{v}\right] =\left[ X,Y\right] ^{v},%
\left[ X^{c},Y^{c}\right] =\left[ X,Y\right] ^{c},%
\end{array}
\\
iii) &
\begin{array}{l}
\chi (U)=Sp\left\{ \frac{\partial }{\partial z^{\alpha }},\frac{\partial }{%
\partial \overline{z}^{\alpha }}\right\} ,\,\chi (TU)=Sp\left\{ \frac{%
\partial }{\partial z^{\alpha }},\frac{\partial }{\partial \overline{z}%
^{\alpha }},\frac{\partial }{\partial z^{\acute{}\alpha }},\frac{\partial }{%
\partial \overline{z}^{\acute{}\alpha }}\right\} , \\
\,(\frac{\partial }{\partial z^{\alpha }})^{c}=\frac{\partial }{\partial
z^{\alpha }},(\frac{\partial }{\partial \overline{z}^{\alpha }})^{c}=\frac{%
\partial }{\partial \overline{z}^{\alpha }},%
\end{array}
\\
iv) &
\begin{array}{l}
\;(\omega +\theta )^{c}=\omega ^{c}+\theta ^{c},\omega ^{c}(Z^{c})=(\omega
Z)^{c}, \\
\omega ^{c}(Z^{v})=\omega ^{v}(Z^{c})=(\omega Z)^{v},\omega
^{c}(Z^{c})=(\omega Z)^{c},%
\end{array}
\\
v) &
\begin{array}{l}
\,\chi ^{*}(U)=Sp\left\{ dz^{\alpha },d\overline{z}^{\alpha }\right\} ,\chi
^{*}(TU)=Sp\left\{ dz^{\alpha },d\overline{z}^{\alpha },dz^{\acute{}\alpha
},d\overline{z}^{\acute{}\alpha }\right\} , \\
(dz^{\alpha })^{c}=\overline{d}z^{\acute{}\alpha },(d\overline{z}^{\alpha
})^{c}=\overline{d}\overline{z}^{\acute{}\alpha },%
\end{array}
\\
vi) & F^{v}(Z^{c})=F^{c}(Z^{v})=(F(Z))^{v},F^{c}(Z^{c})=(F(Z))^{c},%
\end{array}
\]
$\,$

for all\thinspace \thinspace \thinspace $f,g\in \mathcal{F}(M)\,,$ $X,Y,Z\in
\chi (M)$, $\omega ,\theta \in \chi ^{*}(M),\,F\in \Im _{1}^{1}(M).$ Where $%
1\leq \alpha \leq m$ , $\overline{d}\,\,$denotes the differential operator on%
$\,\,TM$ and $\left[ ,\right] $ is Lie bracket.

\subsection{Horizontal Lifts}

The \textit{horizontal\thinspace \thinspace lift} of\thinspace \thinspace $%
f\in \mathcal{F}(M)$ \thinspace to $TM$ is the function $f^{h}\in \mathcal{F}%
(TM)$ given by
\begin{equation}
f^{h}=f^{c}-\gamma (\nabla f),\,\,\,\,\,(\gamma (\nabla f)=\nabla _{\gamma
}f),  \label{2.19}
\end{equation}

where $\nabla $ is an affine linear connection on $M$ with local components $%
\Gamma _{\alpha }^{\beta },$ $\nabla f$ is gradient of $f$ and $\gamma $ is
an operator given by
\begin{equation}
\gamma :\Im _{s}^{r}(M)\rightarrow \Im _{s-1}^{r}(TM).  \label{2.20}
\end{equation}

Thus, it is $f^{h}=0$ since
\begin{equation}
\nabla _{\gamma }f=z^{\acute{}\alpha }(\frac{\partial f}{\partial z^{\alpha }%
})^{v}+\overline{z}^{\acute{}\alpha }(\frac{\partial f}{\partial \overline{z}%
^{\alpha }})^{v}.  \label{2.21}
\end{equation}

The \textit{horizontal lift} of a vector field\textbf{\ }$Z\in \chi (M)$ to $%
TM$ is the vector field $Z^{h}\in \chi (TM)$ given by
\begin{equation}
Z^{h}f^{v}=(Zf)^{v}.  \label{2.22}
\end{equation}

Obviously, we have

\begin{equation}
Z^{h}=Z^{\alpha }D_{\alpha }+\overline{Z}^{\alpha }\overline{D}_{\alpha },{\
}1\leq \alpha ,\beta \leq m  \label{2.23}
\end{equation}

such that $D_{\alpha }=\frac{\partial }{\partial z^{\alpha }}-\Gamma _{\beta
}^{\alpha }\frac{\partial }{\partial z^{\acute{}\alpha }}$ and $\overline{D}%
_{\alpha }=\frac{\partial }{\partial \overline{z}^{\alpha }}-\overline{%
\Gamma }_{\beta }^{\alpha }\frac{\partial }{\partial \overline{z}^{\acute{}%
\alpha }},1\leq \alpha ,\beta \leq m$, where $Z=Z^{\alpha }\frac{\partial }{%
\partial z^{\alpha }}+\overline{Z}^{\alpha }\frac{\partial }{\partial
\overline{z}^{\alpha }}.$

If $\omega =\omega _{\alpha }dz^{\alpha }+\overline{\omega }_{\alpha }d%
\overline{z}^{\alpha }$ we obtain

\begin{equation}
\omega ^{h}=\omega _{\alpha }\eta ^{\alpha }+\overline{\omega }_{\alpha }%
\overline{\eta }^{\alpha },1\leq \alpha ,\beta \leq m  \label{2.24}
\end{equation}

such that $\eta ^{\alpha }=\overline{d}z^{\acute{}\alpha }+\Gamma _{\beta
}^{\alpha }\overline{d}z^{\alpha },$ $\overline{\eta }^{\alpha }=\overline{d}%
\overline{z}^{\acute{}\alpha }+\overline{\Gamma }_{\beta }^{\alpha }%
\overline{d}\overline{z}^{\alpha }.$

The \textit{horizontal lift} of a paracomplex tensor field of type (1,1)%
\textbf{\ } $F\in \Im _{1}^{1}(M)$ to $TM$ is the tensor field $F^{h}\in \Im
_{1}^{1}(TM)$ given by
\begin{equation}
F^{h}(Z^{h})=(FZ)^{h},F^{h}(Z^{v})=(F(Z))^{v}.  \label{2.25}
\end{equation}

The horizontal lift of the paracomplex tensor field of type (1,1)\textbf{\ }$%
F$ is
\begin{equation}
\begin{array}{l}
F^{h}=F_{\alpha }^{\beta }\frac{\partial }{\partial z^{\beta }}\otimes
dz^{\alpha }+(\Gamma _{\beta }^{\alpha }F_{\alpha }^{\beta }-\Gamma _{\alpha
}^{\beta }F_{\beta }^{\alpha })\frac{\partial }{\partial z^{\acute{}\beta }}%
\otimes dz^{\alpha }+F_{\alpha }^{\beta }\frac{\partial }{\partial z^{\acute{%
}\beta }}\otimes dz^{\acute{}\alpha } \\
\,\,\,\,\,\,\,\,\,\,\,+\overline{F}_{\alpha }^{\beta }\frac{\partial }{%
\partial \overline{z}^{\beta }}\otimes d\overline{z}^{\alpha }+(\overline{%
\Gamma }_{\beta }^{\alpha }\overline{F}_{\alpha }^{\beta }-\overline{\Gamma }%
_{\alpha }^{\beta }\overline{F}_{\beta }^{\alpha })\frac{\partial }{\partial
\overline{z}^{\acute{}\beta }}\otimes d\overline{z}^{\alpha }+\overline{F}%
_{\alpha }^{\beta }\frac{\partial }{\partial \overline{z}^{\acute{}\beta }}%
\otimes d\overline{z}^{\acute{}\alpha },%
\end{array}
\label{2.26}
\end{equation}

where $1\leq \alpha ,\beta \leq m,\,F=F_{\alpha }^{\beta }\frac{\partial }{%
\partial z^{\beta }}\otimes dz^{\alpha }+\overline{F}_{\alpha }^{\beta }%
\frac{\partial }{\partial \overline{z}^{\beta }}\otimes d\overline{z}%
^{\alpha }$

With respect to the adapted frame, we have

\begin{equation}
F^{h}=F_{\alpha }^{\beta }D_{\beta }\otimes \theta ^{\alpha }+F_{\alpha
}^{\beta }V_{\beta }\otimes \eta ^{\alpha }+\overline{F}_{\alpha }^{\beta }%
\overline{D}_{\beta }\otimes \overline{\theta }^{\alpha }+\overline{F}%
_{\alpha }^{\beta }\overline{V}_{\beta }\otimes \overline{\eta }^{\alpha },
\label{2.27}
\end{equation}
The properties of horizontal lifts of paracomplex tensor fields are

\[
\begin{array}{ll}
i) & (f+g)^{h}=0,\,\,(f.g)^{h}=0, \\
ii) &
\begin{array}{l}
\,(Z+W)^{h}=Z^{h}+W^{h},\,Z^{h}(f^{v})=(Zf)^{v},\, \\
(\frac{\partial }{\partial z^{\alpha }})^{h}=D_{\alpha },\,(\frac{\partial }{%
\partial \overline{z}^{\alpha }})^{h}=\overline{D}_{\alpha },%
\end{array}
\\
iii) &
\begin{array}{l}
(\omega +\theta )^{h}=\omega ^{h}+\theta ^{h},\,\omega ^{h}(Z^{h})=0,{\ }%
\omega ^{h}(Z^{v})=(\omega Z)^{v},\, \\
(dz^{\alpha })^{h}=\eta ^{\alpha },{\ }(d\overline{z}^{\alpha })^{h}=%
\overline{\eta }^{\alpha },%
\end{array}
\\
iv) & F^{h}(Z^{h})=(FZ)^{h},\,F^{h}(Z^{v})=(F(Z))^{v},%
\end{array}
\]
for all\thinspace \thinspace \thinspace $f,g\in \mathcal{F}(M),$ $Z,W\in
\chi (M),\,\omega ,\theta \in \chi ^{*}(M),\,F\in \Im _{1}^{1}(M).$

Where $1\leq \alpha \leq m,\,\overline{d}\,\,$denotes the differential
operator on$\,\,TM$ and $\,\left[ ,\right] $ is Lie bracket. The set of
local vector fields $\left\{ D_{\alpha },\overline{D}_{\alpha },V_{\alpha }=%
\frac{\partial }{\partial z^{\acute{}\alpha }}\overline{V}_{\alpha }= \frac{%
\partial }{\partial \overline{z}^{\acute{}\alpha }}\right\} $ is called
\textit{adapted frame }to $\nabla .$ The dual coframe $\left\{ \theta
^{\alpha }=dz^{\alpha },\overline{\theta }^{\alpha }=d\overline{z}^{\alpha
},\eta ^{\alpha },\overline{\eta }^{\alpha }\right\} $ is called \textit{%
adapted coframe }to $\nabla .$

\section{\textbf{Lifts of Almost Lorentzian r-Paracontact Structures}}

In this section, we produce almost paracomplex structures on tangent bundle $%
T(\overline{M})$ of almost Lorentzian $r$- paracontact manifold $\overline{M}
$ having the structure $({\mathbf{F}},\xi _{\alpha },\eta ^{\alpha }).$

\subsection{Complete Lifts}

\textbf{Theorem 4.1.} Let $\overline{M}$ be a differentiable manifold
endowed with almost r-paracontact structure $({\mathbf{F}},\xi _{\alpha
},\eta ^{\alpha })$, then
\[
\widetilde{J}={\mathbf{F}}^{c}+\sum_{\alpha =1}^{r}\xi _{\alpha }^{v}\otimes
\eta ^{\alpha v}-\xi _{\alpha }^{c}\otimes \eta ^{\alpha c}
\]

is almost paracomplex structure on $T(\overline{M})$.

\textbf{Proof :} From (\ref{1.10}) and the vertical and complete lifts of
paracomplex tensor fields we have
\begin{equation}
({\mathbf{F}}^{2})^{c}=(I+\sum_{\alpha =1}^{r}\xi _{\alpha }\otimes \eta
_{\alpha })^{c},  \label{3.1}
\end{equation}

\begin{equation}
({\mathbf{F}}^{c})^{2}=I+\sum_{\alpha =1}^{r}\xi _{\alpha }^{v}\otimes \eta
^{\alpha c}+\xi _{\alpha }^{c}\otimes \eta ^{\alpha v},  \label{3.2}
\end{equation}

and
\begin{equation}
{\mathbf{F}}^{c}(\xi _{\alpha }^{v})=0,{\mathbf{F}}^{c}(\xi _{\alpha
}^{c})=0,  \label{3.3}
\end{equation}
\begin{equation}
\eta ^{\alpha v}\circ {\mathbf{F}}^{c}=0,\eta ^{\alpha c}\circ {\mathbf{F}}%
^{v}=0,\eta ^{\alpha c}\circ {\mathbf{F}}^{c}=0,  \label{3.4}
\end{equation}
\begin{equation}
\eta ^{\alpha v}(\xi _{\beta }^{v})=0,\eta ^{\alpha v}(\xi _{\beta
}^{c})=\delta _{\beta }^{\alpha },\eta ^{\alpha c}(\xi _{\beta }^{v})=\delta
_{\beta }^{\alpha },\,\eta ^{\alpha c}(\xi _{\beta }^{c})=0.  \label{3.5}
\end{equation}

Consider an element $\widetilde{J}$ of $\Im _{1}^{1}(T\overline{M})$ given by

\begin{equation}
\widetilde{J}={\mathbf{F}}^{c}+\sum_{\alpha =1}^{r}(\xi _{\alpha
}^{v}\otimes \eta ^{\alpha v}-\xi _{\alpha }^{c}\otimes \eta ^{\alpha c}).
\label{3.6}
\end{equation}

Using (\ref{3.2}) and (\ref{3.6}) we find equation
\begin{eqnarray}
(\tilde{J})^{2} &=&({\mathbf{F}}^{c}+\sum_{\alpha =1}^{r}\xi _{\alpha
}^{v}\otimes \eta ^{\alpha v}-\xi _{\alpha }^{c}\otimes \eta ^{\alpha c})^{2}
\label{3.7} \\
&=&({\mathbf{F}}^{c})^{2}+\sum_{\alpha =1}^{r}\,{\mathbf{F}}^{c}(\xi
_{\alpha }^{v}\otimes \eta ^{\alpha v}-\xi _{\alpha }^{c}\otimes \eta
^{\alpha c})+\sum_{\alpha =1}^{r}\,(\xi _{\alpha }^{v}\otimes \eta ^{\alpha
v}-\xi _{\alpha }^{c}\otimes \eta ^{\alpha c}){\mathbf{F}}^{c}  \nonumber \\
&&+\sum_{\alpha =1}^{r}(\xi _{\alpha }^{v}\otimes \eta ^{\alpha v}-\xi
_{\alpha }^{c}\otimes \eta ^{\alpha c})^{2}  \nonumber \\
&=&I+\sum_{\alpha =1}^{r}(\xi _{\alpha }^{v}\otimes \eta ^{\alpha c}+\xi
_{\alpha }^{c}\otimes \eta ^{\alpha v})+  \nonumber \\
&&+\sum_{\alpha =1}^{r}({\mathbf{F}}^{c}(\xi _{\alpha }^{v})\eta ^{\alpha v}-%
{\mathbf{F}}^{c}(\xi _{\alpha }^{c})\eta ^{\alpha c}+(\eta ^{\alpha v}\circ {%
\mathbf{F}}^{c})\xi _{\alpha }^{v}-(\eta ^{\alpha c}\circ {\mathbf{F}}%
^{c})\xi _{\alpha }^{c}  \nonumber \\
&&+\sum_{\alpha =1}^{r}(\xi _{\alpha }^{v}\otimes (\eta ^{\alpha v}(\xi
_{\alpha }^{v}))\eta ^{\alpha v}-\xi _{\alpha }^{v}\otimes (\eta ^{\alpha
v}(\xi _{\alpha }^{c}))\eta ^{\alpha c}  \nonumber \\
&&-\xi _{\alpha }^{c}\otimes (\eta ^{\alpha c}(\xi _{\alpha }^{v}))\eta
^{\alpha v}+\xi _{\alpha }^{c}\otimes (\eta ^{\alpha c}(\xi _{\alpha
}^{c}))\eta ^{\alpha c}).  \nonumber
\end{eqnarray}
By means of (\ref{3.3}), (\ref{3.4}) and (\ref{3.5}), we obtain

\begin{equation}
(\widetilde{J})^{2}=I.  \label{3.8}
\end{equation}
So, $\widetilde{J}$ is an almost paracomplex structure in $T(\bar{M})$.
Hence the proof is completed.

\textbf{Theorem 4.2.} Let $\overline{M}$ be a differentiable manifold
endowed with almost Lorentzian r-paracontact structure $({\mathbf{F}},\xi
_{\alpha },\eta ^{\alpha })$, then an almost paracomplex structure on $T(%
\overline{M})$ is calculated by
\[
\widehat{J}={\mathbf{F}}^{c}-\sum_{\alpha =1}^{r}\xi _{\alpha }^{v}\otimes
\eta ^{\alpha v}-\xi _{\alpha }^{c}\otimes \eta ^{\alpha c}.
\]

\textbf{Proof :} By means of the equation (\ref{1.14}) and the vertical and
complete lifts of paracomplex tensor fields we have
\begin{equation}
({\mathbf{F}}^{2})^{c}=(I-\sum_{\alpha =1}^{r}\xi _{\alpha }\otimes \eta
_{\alpha })^{c},  \label{3.9}
\end{equation}

\begin{equation}
({\mathbf{F}}^{c})^{2}=I-\sum_{\alpha =1}^{r}\xi _{\alpha }^{v}\otimes \eta
^{\alpha c}+\xi _{\alpha }^{c}\otimes \eta ^{\alpha v},  \label{3.10}
\end{equation}

and
\begin{equation}
{\mathbf{F}}^{c}(\xi _{\alpha }^{v})=0,{\mathbf{F}}^{c}(\xi _{\alpha
}^{c})=0,  \label{3.11}
\end{equation}
\begin{equation}
\eta ^{\alpha v}\circ {\mathbf{F}}^{c}=0,\eta ^{\alpha c}\circ {\mathbf{F}}%
^{v}=0,\eta ^{\alpha c}\circ {\mathbf{F}}^{c}=0,  \label{3.12}
\end{equation}
\begin{equation}
\eta ^{\alpha v}(\xi _{\beta }^{v})=0,\eta ^{\alpha v}(\xi _{\beta
}^{c})=-\delta _{\beta }^{\alpha },\eta ^{\alpha c}(\xi _{\beta
}^{v})=-\delta _{\beta }^{\alpha },\,\eta ^{\alpha c}(\xi _{\beta }^{c})=0.
\label{3.13}
\end{equation}

Take an element $\widehat{J}$ of $\Im _{1}^{1}(T\overline{M})\;$defined by

\begin{equation}
\widehat{J}={\mathbf{F}}^{c}-\sum_{\alpha =1}^{r}(\xi _{\alpha }^{v}\otimes
\eta ^{\alpha v}-\xi _{\alpha }^{c}\otimes \eta ^{\alpha c})  \label{3.15}
\end{equation}

Similarly proof of \textbf{Theorem 4.1,} using (\ref{3.2}), (\ref{3.3}), (%
\ref{3.4}) and (\ref{3.5})and (\ref{3.6}) we have the equation

\begin{equation}
(\widehat{J})^{2}=I.  \label{3.16}
\end{equation}
Thus, for $\widehat{J}$ is an almost paracomplex structure in $T(\bar{M})$,
the proof is completed.

\subsection{Horizontal Lifts}

\textbf{Theorem 4.3.} Let $({\mathbf{F}},\xi _{\alpha },\eta ^{\alpha })$ be
an almost r-paracontact structure in $\overline{M}$ with an affine
connection $\nabla $. Then an almost paracomplex in $T(\overline{M})$ is
given by
\[
\widetilde{J}^{*}={\mathbf{F}}^{h}+\sum_{\alpha =1}^{r}\xi _{\alpha
}^{v}\otimes \eta ^{\alpha v}-\xi _{\alpha }^{h}\otimes \eta ^{\alpha h}).
\]
\textbf{Proof :} Taking into consideration the equation given by (\ref{1.10}%
) and the horizontal lifts of paracomplex tensor fields, we have
\begin{equation}
({\mathbf{F}}^{2})^{h}=(I+\sum_{\alpha =1}^{r}\xi _{\alpha }\otimes \eta
^{\alpha })^{h},  \label{3.17}
\end{equation}
\begin{equation}
({\mathbf{F}}^{h})^{2}=I+\sum_{\alpha =1}^{r}\xi _{\alpha }^{h}\otimes \eta
^{\alpha v}+\xi _{\alpha }^{v}\otimes \eta ^{\alpha h},  \label{3.18}
\end{equation}
and
\begin{equation}
{\mathbf{F}}^{h}(\xi _{\alpha }^{h})=0,{\mathbf{F}}^{h}(\xi _{\alpha
}^{c})=0,  \label{3.19}
\end{equation}
\begin{equation}
\eta ^{\alpha h}\circ {\mathbf{F}}^{h}=0,\eta ^{\alpha v}\circ {\mathbf{F}}%
^{h}=0.  \label{3.20}
\end{equation}
\begin{equation}
\eta ^{\alpha h}(\xi _{\beta }^{h})=0,\eta ^{\alpha h}(\xi _{\beta
}^{v})=\delta _{\beta }^{\alpha },\eta ^{\alpha v}(\xi _{\beta }^{h})=\delta
_{\beta }^{\alpha }.\,  \label{3.21}
\end{equation}
Given an element $\widetilde{J}$ of $\Im _{1}^{1}(T\overline{M})\;$defined
by
\begin{equation}
\widetilde{J}^{*}={\mathbf{F}}^{h}+{\sum_{\alpha =1}^{r}}(\xi _{\alpha
}^{v}\otimes \eta ^{\alpha v}-\xi _{\alpha }^{h}\otimes \eta ^{\alpha h}).
\label{3.22}
\end{equation}
Taking care of the above equations, it is clear that
\begin{equation}
(\widetilde{J}^{*})^{2}=I.  \label{3.23}
\end{equation}
Consequently $\widetilde{J}^{*}$ is an almost paracomplex structure in $T(%
\overline{M})$. Thus the proof is finished. \thinspace \thinspace \thinspace
\thinspace \thinspace \thinspace \thinspace \thinspace \thinspace \thinspace

\textbf{Theorem 4.4.} Let $({\mathbf{F}},\xi _{\alpha },\eta ^{\alpha })$ be
an almost Lorentzian r-paracontact structure in $\overline{M}$ with an
affine connection $\nabla $. Then structure
\[
\widehat{J}^{*}={\mathbf{F}}^{h}-\sum_{\alpha =1}^{r}\xi _{\alpha
}^{v}\otimes \eta ^{\alpha v}-\xi _{\alpha }^{h}\otimes \eta ^{\alpha h}).
\]

is an almost paracomplex in $T(\overline{M}).$

\textbf{Proof:} It is similar to the proofs of the above theorems.

\section{\textbf{Corollary}}

Taking into consideration the above theorems, we conclude that when we
consider an almost Lorentzian $r$- paracontact structure on the base
manifold, the structure defined on the tangent bundle $T(\overline{M})$ is
an almost paracomplex. Therefore, by means of [13], it is possible to obtain
paracomplex Hamiltonian formalisms in classical mechanics and field theory
on $T(\overline{M}),$ is the tangent manifold of a differentiable manifold $%
\overline{M}$ endowed with almost Lorentzian r-paracontact structure $({%
\mathbf{F}},\xi _{\alpha },\eta ^{\alpha })$.

\textbf{REFERENCES}

[1] M. Tekkoyun, \c{S}. Civelek, First order lifts of complex structures,
Algebras Groups and Geometries (AGG), \textbf{19}(3), (2002)373-382.

[2] M. Tekkoyun, \c{S}. Civelek, On lifts of structures on complex
manifolds, Differential Geometry-Dynamics Systems,\textbf{5}(1) (2003)59-64.

[3] M. Tekkoyun, On horizontal lifts of complex structures, Hadronic Journal
Supplement (HJS), \textbf{18}(4)(2003)411-424

[4] M. Tekkoyun, On lifts of paracomplex structures, Turk. J. Math.,\textbf{%
30}(2006)197-210.

[5] P. Dombrowski, On the geometry of the tangent bundles, Jour. Reine und
Angew. Math., \textbf{210} (1962)73-88.

[6] S. Ianus, C. Udriste, On the tangent bundle of a differentiable
manifold, Stud. Cerc. Mat., \textbf{22} (1970), 599-611.

[7] K. Yano, S. Ishihara, Differential geometry in tangent bundles, Kodai
Math. Sem. Rep.,\textbf{18}(1996)271-292.

[8] D.E. Blair, Contact manifolds in Riemannian geometry, Lecture Notes in
Math, 509, Springer Verlag, New York, (1976).

[9] L.S. Das, M.N.I. Khan, Almost r-contact structures on the tangent
bundle, Differential Geometry-Dynamics Systems, \textbf{7}(2005)34-41.

[10] K. Matsumoto, On Lorentzian paracontact manifolds, Bull. of Yamagata
Univ. Nat. Sci, \textbf{12}(2) (1989)151-156.

[11] A. J. Ledger, K. Yano, Almost complex structures on complex bundles,
Jour. Diff. Geometry, \textbf{1}(1967), 355-368.

[12] K. Yano, S. Ishihara, Almost complex structures induced in tangent
bundles, Kodai Math. Sem. Rep., \textbf{19}(1967)1-27.

[13] M. Tekkoyun, On Para-Euler Lagrange and para- Hamiltonian equations,
Physics Letters A, \textbf{340}, (2005)7-11.

\end{document}